\documentclass[11pt,letterpaper]{article}
\usepackage{amsfonts,amsmath,calc,graphicx}
\allowdisplaybreaks
\usepackage[sort&compress,numbers]{natbib}

\DeclareFontShape{OT1}{cmr}{b}{sc}
{<5><6><7><8><9><10><12><10.95><14.4><17.28><20.74><24.88>cmbcsc10}{}
\DeclareFontShape{OT1}{cmr}{bx}{sc}
{<5><6><7><8><9><10><12><10.95><14.4><17.28><20.74><24.88>cmbcsc10}{}

\DeclareFontShape{OT1}{cmr}{b}{tt}
{<5><6><7><8><9><10><12><10.95><14.4><17.28><20.74><24.88>cmbtt10}{}
\DeclareFontShape{OT1}{cmr}{bx}{tt}
{<5><6><7><8><9><10><12><10.95><14.4><17.28><20.74><24.88>cmbtt10}{}

\newcommand{\FancyFootnotes}{\renewcommand{\thefootnote}{\fnsymbol{footnote}}}
\newcommand{\FancyFootnotesOff}{\renewcommand{\thefootnote}{\arabic{footnote}}}

\newcommand{\paran}[1]{\textup{(}#1\textup{)}}

\newcommand{\bracket}[1]{\ensuremath{\protect\left(#1\right)}}

\newcommand{\ceil}[1]{\ensuremath{\protect\lceil#1\rceil}}

\newcommand{\half}{\ensuremath{\protect\tfrac{1}{2}}}

\newcommand{\threehalves}{\ensuremath{\protect\tfrac{3}{2}}}

\newcommand{\quarter}{\ensuremath{\protect\tfrac{1}{4}}}

\newcommand{\third}{\ensuremath{\protect\tfrac{1}{3}}}

\newcommand{\Oh}[1]{\ensuremath{\protect\mathcal{O}(#1)}}

\newcommand{\Figure}[4][htb]{
\begin{figure}[#1]
	\vspace*{1ex}
	\begin{center}#3\end{center}
	\vspace*{-3ex}
	\caption{\figlabel{#2}#4}
\end{figure}
}

\newcommand{\Settings}[5]{
\renewcommand{\baselinestretch}{#1}
\setlength{\footnotesep}{\baselinestretch\footnotesep}

\setlength{\marginparwidth}{0mm}
\setlength{\marginparsep}{0mm}
\setlength{\marginparpush}{0mm}

\newlength{\LRmargin}
\setlength{\LRmargin}{#2}
\newlength{\TBmargin}
\setlength{\TBmargin}{#3}

\setlength{\voffset}{-1in}           % cancels default printer margin
\setlength{\hoffset}{-1in}           % cancels default printer margin

\setlength{\headsep}{#4}
\setlength{\footskip}{#5}

\setlength{\oddsidemargin}{\LRmargin}
\setlength{\evensidemargin}{\LRmargin}
\setlength{\topmargin}{\TBmargin}

\setlength{\textheight}{\paperheight}
\addtolength{\textheight}{-1in}			% default printer margin
\addtolength{\textheight}{-1.0\voffset}
\addtolength{\textheight}{-2.0\TBmargin}	% top and bottom
\addtolength{\textheight}{-1.0\headheight}
\addtolength{\textheight}{-1.0\headsep}
\addtolength{\textheight}{-1.0\footskip}

\setlength{\textwidth}{\paperwidth}
\addtolength{\textwidth}{-1.0in}		% default printer margin
\addtolength{\textwidth}{-1.0\hoffset}
\addtolength{\textwidth}{-2.0\LRmargin}		% left and right margin
}

\newlength{\marginboxwidth}

\usepackage{pifont}

\usepackage{amsthm}

\theoremstyle{plain}

\newtheorem{theorem}{Theorem}

\newtheorem{lemma}{Lemma}
\newtheorem{corollary}{Corollary}

\theoremstyle{definition}
\newtheorem{conjecture}{Conjecture}

\newtheorem{open}{Open Problem}

\newcommand{\twothmref}[2]{Theorems~\ref{thm:#1} and \ref{thm:#2}}
\newcommand{\lemref}[1]{Lemma~\ref{lem:#1}}

\newcommand{\thmref}[1]{Theorem~\ref{thm:#1}}

\newcommand{\figref}[1]{Figure~\ref{fig:#1}}
\newcommand{\secref}[1]{Section~\ref{sec:#1}}
\newcommand{\twosecref}[2]{Sections~\ref{sec:#1} and \ref{sec:#2}}

\newcommand{\corref}[1]{Corollary~\ref{cor:#1}}

\newcommand{\eqnref}[1]{\eqref{eqn:#1}}

\newcommand{\seclabel}[1]{\label{sec:#1}}
\newcommand{\figlabel}[1]{\label{fig:#1}}

\newcommand{\thmlabel}[1]{\label{thm:#1}}
\newcommand{\lemlabel}[1]{\label{lem:#1}}

\newcommand{\corlabel}[1]{\label{cor:#1}}
\newcommand{\eqnlabel}[1]{\label{eqn:#1}}

\newcommand{\mySection}[2]{\section{#1}\seclabel{#2}}

\usepackage{charter}

\newcommand{\OP}{\ensuremath{\mathcal{OP}}}
\newcommand{\G}{\ensuremath{\mathcal{G}}}

\Settings{1.2}{30mm}{25mm}{0cm}{1cm}

\FancyFootnotes

\begin{document}

\title{\textbf{Induced Subgraphs of Bounded Degree\\ and Bounded Treewidth}\footnote{A short version of this paper will appear in \emph{Proc.\ of  31st International Workshop on Graph Theoretic Concepts in Computer Science} (WG '05), Lecture Notes in Computer Science, Springer.}}

\author{Prosenjit Bose\,\footnotemark[4] 
\and Vida Dujmovi{\'c}\,\footnotemark[4] 
\and David R. Wood\,\footnotemark[5]}

\footnotetext[4]{School of Computer Science, Carleton University, Ottawa, Canada (\texttt{\{jit,vida\}@scs.carleton.ca}). Research supported by NSERC.}

\footnotetext[5]{Departament de Matem{\`a}tica Aplicada II, Universitat Polit{\`e}cnica de Catalunya, Barcelona, Spain (\texttt{david.wood@upc.edu}).
Supported by Government of Spain grant MEC SB2003-0270. Research partially completed at Carleton Univerity.}

\maketitle

\begin{abstract}
We prove that for all $0\leq t\leq k$ and $d\geq 2k$, every graph $G$ with treewidth at most $k$ has a `large' induced subgraph $H$, where $H$ has treewidth at most $t$ and every vertex in $H$ has degree at most $d$ in $G$. 
The order of $H$ depends on $t$, $k$, $d$, and the order of $G$.
With $t=k$, we obtain large sets of bounded degree vertices. 
With $t=0$, we obtain large independent sets of bounded degree. In both these cases, our bounds on the order of $H$ are tight. For bounded degree independent sets in trees, we characterise the extremal graphs. Finally, we prove that an interval graph with maximum clique size $k$ has a maximum independent set in which every vertex has degree at most $2k$.
\end{abstract}

\FancyFootnotesOff

%%%%%%%%%%%%%%%%%%%%%%%%%%%%%%%%%%%%%%%%%%%%%%%%%%%%%%%%%%%%%%%%%%%%%%%%%%%%
\section{Introduction}\seclabel{Intro}
%%%%%%%%%%%%%%%%%%%%%%%%%%%%%%%%%%%%%%%%%%%%%%%%%%%%%%%%%%%%%%%%%%%%%%%%%%%%

The `treewidth' of a graph has arisen as an important parameter in the Robertson/Seymour theory of graph minors and in algorithmic complexity. See \citet{Bodlaender-TCS98} and \citet{Reed-AlgoTreeWidth03} for surveys on treewidth. The main result of this paper, proved in \secref{LargeBoundedTreewidthDegreeSubgraphs}, states that every graph $G$ has a large induced subgraph of bounded treewidth in which every vertex has bounded degree in $G$. The order of the subgraph depends on the treewidth of $G$, the desired treewidth of the subgraph, and the desired degree bound. Moreover, we prove that the bound is best possible in a number of cases.

Before that, in \twosecref{LargeSubgraphsBoundedDegree}{LargeSubgraphsBoundedTreewidth} we consider two relaxations of the main result, firstly without the treewidth constraint, and then without the degree constraint. That is, we determine the minimum number of vertices of bounded degree in a graph of given treewidth (\secref{LargeSubgraphsBoundedDegree}), and we determine the minimum number of vertices in an induced subgraph of bounded treewidth, taken over all graphs of given treewidth (\secref{LargeSubgraphsBoundedTreewidth}). This latter result is the first ingredient in the proof of the main result. The second ingredient is proved in \secref{Structure}, where we consider the structure of the subgraph of a $k$-tree induced by the vertices of bounded degree. In particular, we prove that this subgraph has surprisingly small treewidth.

A graph with treewidth $0$ has no edges. Thus our results pertain to independent sets for which every vertex has bounded degree in $G$. Here our bounds are tight, and in the case of trees, we characterise the extremal trees. Furthermore, by exploiting some structural properties of interval graphs that are of independent interest, we prove that every interval graph with no $(k+2)$-clique has a maximum independent set in which every vertex has degree at most $2k$. These results are presented in \secref{IndSets}. 

\subsection{Preliminaries}
%%%%%%%%%%%%%%%%%%%%%%%%%%%%%%%%%%%%%%%%%%%%%%%%%%%%%%%%%%%%%%%%%%%%%%%%%%%%

Let $G$ be a graph. All graphs considered are finite, undirected, and simple. 
The vertex-set and edge-set of $G$ are denoted by $V(G)$ and $E(G)$,
respectively. The number of vertices of $G$ is denoted by $n=|V(G)|$. The subgraph \emph{induced} by a set of vertices $S\subseteq V(G)$ has vertex set $S$ and edge set $\{vw\in E(G):v,w\in S\}$, and is denoted by $G[S]$. 

A $k$-\emph{clique} ($k\geq0$) is a set of $k$ pairwise adjacent vertices. Let $\omega(G)$ denote the maximum number $k$ such that $G$ has a $k$-clique. A \emph{chord} of a cycle $C$ is an edge not in $C$ whose endpoints are both in $C$. $G$ is \emph{chordal} if every cycle on at least four vertices has a chord. The \emph{treewidth} of $G$ is the minimum number $k$ such that $G$ is a subgraph of a chordal graph $G'$ with $\omega(G')\leq k+1$. 

A vertex is \emph{simplicial} if its neighbourhood is a clique. For each vertex $v\in V(G)$, let $G\setminus v$ denote the subgraph $G[V(G)\setminus\{v\}]$. The family of graphs called $k$-trees $(k\geq0$) are defined recursively as follows. A graph $G$ is a $k$-\emph{tree} if one of the following conditions are satisfied:
\begin{enumerate}
\item[(a)] $G$ is a ($k+1$)-clique, or
\item[(b)] $G$ has a simplicial vertex $v$ whose neighbourhood is a
$k$-clique, and $G\setminus v$ is a $k$-tree.
\end{enumerate}

By definition, the graph obtained from a $k$-tree $G$ by adding a new vertex $v$ adjacent to each vertex of a $k$-clique $C$ is also a $k$-tree, in which case we say $v$ is \emph{added onto} $C$. Every $k$-tree $G$ on $n$ vertices satisfies the following obvious facts:

\begin{itemize} 

\item $\omega(G)=k+1$

\item $G$ has minimum degree $k$

\item $G$ has $kn-\half k(k+1)$ edges, and thus has average degree $2k-k(k+1)/n$.

\end{itemize}

It is well known that the treewidth of a graph $G$ equals the minimum number $k$ such that $G$ is a spanning subgraph of a $k$-tree. 

We will express our results using the following notation. Let $G$ be a graph. Let $V_d(G)=\{v\in V(G):\deg_G(v)\leq d\}$ denote the set of vertices of $G$ with degree at most $d$. Let $G_d=G[V_d(G)]$. A subset of $V_d(G)$ is called a \emph{degree-$d$} set. For an integer $t\geq0$, a \emph{$t$-set} of $G$ is a set $S$ of vertices of $G$ such that the induced subgraph $G[S]$ has treewidth at most $t$. Let $\alpha^t(G)$ be the maximum number of vertices in a $t$-set of $G$. Let $\alpha^t_d(G)$ be the maximum number of vertices in a degree-$d$ $t$-set of $G$. Observe that $\alpha^t_d(G)=\alpha^t(G_d)$. 

Let $\G$ be a family of graphs. Let $\alpha^t(\G)$ be the minimum of $\alpha^t(G)$, and let $\alpha^t_d(\G)$ be the minimum of $\alpha^t_d(G)$, taken over all $G\in\G$. Let $\G_{n,k}$ be the family of $n$-vertex graphs with treewidth $k$. Note that every graph in $\G_{n,k}$ has at least $k+1$ vertices. These definitions imply the following. Every graph $G\in\G$ has
$\alpha_d^t(G)\geq\alpha_d^t(\G)$ and $\alpha^t(G)\geq\alpha^t(\G)$. Furthermore, there is at least one graph $G$ for which $\alpha_d^t(G)=\alpha_d^t(\G)$, and  there is at least one graph $G$ for which $\alpha^t(G)=\alpha^t(\mathcal{G})$. Thus the lower bounds we derive in this paper are universal and the upper bounds are existential.

As described above, our main result is a lower bound on $\alpha^t_d(\G_{n,k})$ that is tight in many cases. Here, lower and upper bounds are 'tight' if they are equal when ignoring the terms independent of $n$. Many of our upper bound constructions are based on the $k$-th power of an $n$-vertex path $P_n^k$. This graph has vertex set $\{v_1,v_2,\dots,v_n\}$ and edge set $\{v_iv_j:|i-j|\leq k\}$. Obviously $P_n^k$ is a $k$-tree.

For $t=k$, a degree-$d$ $t$-set in a graph $G$ with treewidth $k$ is simply a set of vertices with degree at most $d$. Thus in this case, $\alpha_d^k(G)=|V_d(G)|$. At the other extreme, a graph has treewidth $0$ if and only if it has no edges. A set of vertices $I\subseteq V(G)$ is \emph{independent} if $G[I]$ has no edges. Thus a $0$-set of $G$ is simply an independent set of vertices of $G$. As is standard, we abbreviate $\alpha^0(G)$ by $\alpha(G)$, $\alpha_d^0(G)$ by $\alpha_d(G)$, etc\footnote[1]{Be aware that some authors defined a degree-$d$ independent set to consist of vertices of degree \emph{strictly} less than $d$.}. An independent set $I$ of $G$ is \emph{maximum} if $|I|\geq|J|$ for every independent set $J$ of $G$. Thus $\alpha(G)$ is the cardinality of a maximum independent set of $G$. 

%%%%%%%%%%%%%%%%%%%%%%%%%%%%%%%%%%%%%%%%%%%%%%%%%%%%%%%%%%%%%%%%%%%%%%%%%%%%%%%%
\mySection{Large Subgraphs of Bounded Degree}{LargeSubgraphsBoundedDegree}
%%%%%%%%%%%%%%%%%%%%%%%%%%%%%%%%%%%%%%%%%%%%%%%%%%%%%%%%%%%%%%%%%%%%%%%%%%%%%%%%

In this section we prove tight lower bounds on the number of vertices of bounded degree in graphs of treewidth $k$. We will use the following result of \citet{BSW-DM04}.

\begin{lemma}[\citep{BSW-DM04}]
\lemlabel{SizeVd}
Let $G$ be a graph on $n$ vertices, with minimum degree $\delta$, and with average degree $\alpha$. Then for every integer $d\geq\delta$,
\begin{equation*} 
|V_d(G)|\;\geq\;\bracket{\frac{d+1-\alpha}{d+1-\delta}}n\enspace.
\end{equation*}
\end{lemma}

\begin{theorem} 
\thmlabel{kSet}
For all integers $k\geq0$ and $d\geq2k-1$,
\begin{equation*}
\lim_{n\rightarrow\infty}\frac{\alpha^k_d(\G_{n,k})}{n}
\;=\;\frac{d-2k+1}{d-k+1}\enspace.
\end{equation*}
\end{theorem}

\begin{proof} First we prove a lower bound on $\alpha^k_d(\G_{n,k})$. Let $G$ be a graph in $\G_{n,k}$ with $\alpha^k_d(G)=\alpha^k_d(\G_{n,k})$. If a vertex $v$ of $G$ has degree at most $d$ in a spanning supergraph of $G$, then $v$ has degree at most $d$ in $G$. Thus we can assume that $G$ is a $k$-tree. Hence $G$ has minimum degree $k$ and average degree $2k-k(k+1)/n$. By \lemref{SizeVd},
\begin{equation}
\eqnlabel{kSetLowerBound} 
\alpha^k_d(\G_{n,k})
\;=\;
|V_d(G)|
\;\geq\;
\bracket{\frac{d+1-2k+k(k+1)/n}{d+1-k}}n
\;=\;
\bracket{\frac{d-2k+1}{d-k+1}}n\;+\;\frac{k(k+1)}{d-k+1}\enspace.
\end{equation}

Now we prove an upper bound on $\alpha^k_d(\G_{n,k})$ for all $n\equiv2k\pmod{d-k+1}$, and for all $k\geq0$ and $d\geq2k-1$. 
Let $s$ be the integer such that $n-2k=s(d-k+1)$. Then $s\geq0$. We now construct a graph $G\in\G_{n,k}$. Initially let $G=P^k_{(s+2)k}$ be the $k$-th power of the path $(v_1,v_2,\dots,v_{(s+2)k})$. Let $r=d-2k+1$. Then $r\geq0$. Add $r$ vertices onto the clique $(v_{ik+1},v_{ik+2},\dots,v_{ik+k})$ for each $1\leq i\leq s$. Thus $G$ is a $k$-tree, as illustrated in \figref{kSetExample}. The number of vertices in $G$ is 
\begin{equation}
(s+2)k+sr\;=\;(s+2)k+s(d-2k+1)\;=\;s(d-k+1)+2k\;=\;n
\enspace.
\end{equation}
Each vertex $v_i$, $k+1\leq i\leq(s+1)k$, has degree $2k+r=d+1$. Hence such a vertex is not in a degree-$d$ set. The remaining vertices all have degree at most $d$. Thus
\begin{equation}
\eqnlabel{kSetUpperBound} 
\alpha^k_d(\G_{n,k})\;\leq\;
\alpha^k_d(G)\;=\;
|V_d(G)|\;=\;rs+2k%\\
%\;=\;&\frac{(d-2k+1)(n-2k)}{d-k+1}\;+\;2k\\
%\;=\;&\bracket{\frac{d-2k+1}{d-k+1}}n\;+\;
%	2k\bracket{1-\frac{d-2k+1}{d-k+1}}\\
\;=\;\bracket{\frac{d-2k+1}{d-k+1}}n\;+\;
	\frac{2k^2}{d-k+1}
\enspace.
\end{equation}

\Figure[!htb]{kSetExample}{\includegraphics{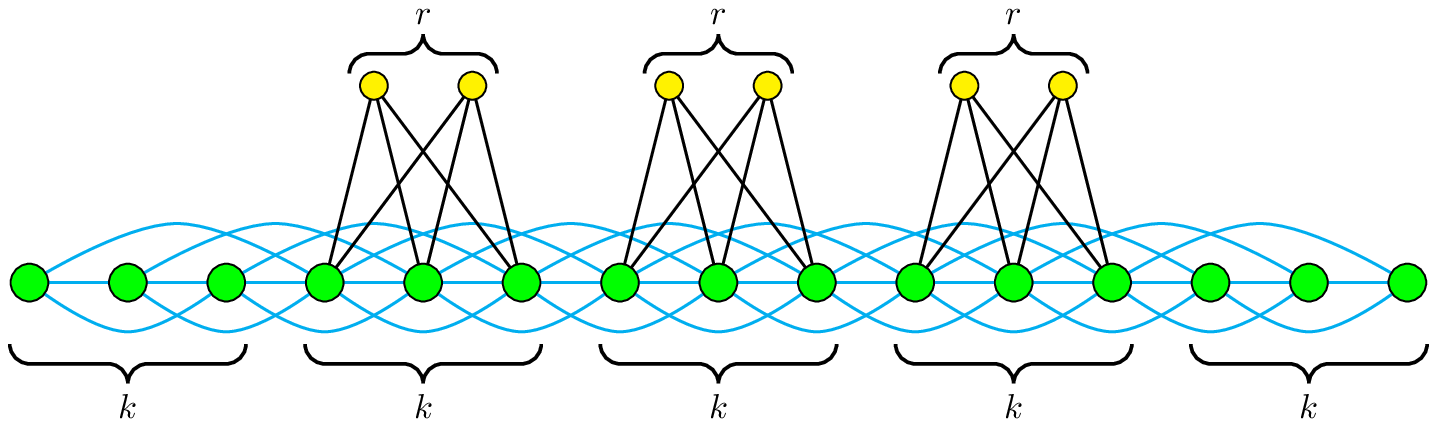}}{The graph $G$ with
$k=3$, $d=7$, and $s=3$ (and thus $r=2$).}

Before we prove the limit that it is claimed in the theorem, note that the difference between the lower and upper bounds in \eqnref{kSetLowerBound} and \eqnref{kSetUpperBound} is only 
\begin{equation*}
\frac{2k^2-k(k+1)}{d+1-k}\;=\;\frac{k(k-1)}{d+1-k}\;\leq\;k-1
\enspace.
\end{equation*}

Given any $n\geq 2k$, there is an integer $n'$ such that $n\leq n'\leq n+d-k$ and $n'\equiv 2k\pmod{d-k+1}$. Hence 
\begin{align*}
\alpha^k_d(\G_{n,k})
\;\leq\;
\alpha^k_d(\G_{n',k})
&\;\leq\;
\bracket{\frac{d-2k+1}{d-k+1}}n'+\frac{2k^2}{d-k+1}\\
&\;\leq\;
\bracket{\frac{d-2k+1}{d-k+1}}n+
\frac{(d-2k+1)(d-k)+2k^2}{d-k+1}
\enspace.
\end{align*}
By \eqnref{kSetLowerBound}, for all $n$,
\begin{align*}
\frac{k(k+1)}{(d-k+1)n}
\;\leq\;
\frac{\alpha^k_d(\G_{n,k})}{n}\,-\,\frac{d-2k+1}{d-k+1}
\;\leq\;
\frac{(d-2k+1)(d-k)+2k^2}{(d-k+1)n}
\enspace.
\end{align*}
Therefore for all $\epsilon>0$, there is an $n_0$ such that for all $n\geq n_0$,
\begin{align*}
0\leq\;\frac{\alpha^k_d(\G_{n,k})}{n}-\frac{d-2k+1}{d-k+1}
\;\leq\;\epsilon\enspace.
\end{align*}
Therefore the sequence $\{\alpha^k_d(\G_{n,k})/n:n\geq 2k\}$ converges to $\frac{d-2k+1}{d-k+1}$.
\end{proof}

%%%%%%%%%%%%%%%%%%%%%%%%%%%%%%%%%%%%%%%%%%%%%%%%%%%%%%%%%%%%%%%%%
\mySection{Large Subgraphs of Bounded Treewidth}{LargeSubgraphsBoundedTreewidth}
%%%%%%%%%%%%%%%%%%%%%%%%%%%%%%%%%%%%%%%%%%%%%%%%%%%%%%%%%%%%%%%%%

We now prove a tight bound on the maximum order of an induced subgraph of bounded treewidth in a graph of treewidth $k$.

\begin{theorem}
\thmlabel{tSet}
For all integers $n$ and $0\leq t\leq k$, 
\begin{equation*}
\alpha^t(\G_{n,k})\;=\;\bracket{\frac{t+1}{k+1}}n\enspace.
\end{equation*}
\end{theorem}

\begin{proof} First we prove the lower bound. Let $G$ be a graph in $\G_{n,k}$. First suppose that $G$ is a $k$-tree. By definition, $V(G)$ can be ordered $(v_1,v_2,\dots,v_n)$ so that for each vertex $v_i$, the \emph{predecessors} $\{v_j:j<i,v_iv_j\in E(G)\}$ of $v_i$ are a clique of $\min\{k,i-1\}$ vertices. Now colour $G$ greedily in this order. That is, for $i=1,2,\dots,n$, assign to $v_i$ the minimum positive integer (a colour) not already assigned to a neighbour of $v_i$. Clearly $k+1$ colours suffice. Let
$S$ be the union of the $t+1$ largest colour classes (monochromatic set of vertices). Thus $|S|\geq(t+1)n/(k+1)$. For each vertex $v_i$ in $S$, the predecessors of $v_i$ that are in $S$ and $v_i$ itself form a clique, and thus have pairwise distinct colours. Thus $v_i$ has at most $t$ predecessors in $S$, and they form a clique in $G[S]$. Hence $G[S]$ has treewidth at most $t$, and $S$ is the desired $t$-set. Now suppose that $G$ is not a $k$-tree. Then $G$ is a spanning subgraph of a $k$-tree $G'$. Thus $G'$ has a $t$-set $S$ with at least $(t+1)n/(k+1)$ vertices. Now $G[S]$ is a subgraph of $G'[S]$. Thus $G[S]$ also has treewidth at most $t$. 

For the upper bound, we now show that every $t$-set of $P_n^k$ has at most $(t+1)n/(k+1)$ vertices. First suppose that $t=0$. A $0$-set is an independent set. Clearly every independent set of $P_n^k$ has at most $n/(k+1)$ vertices. Now consider the case of general $t$. Let $S$ be a $t$-set of $P_n^k$. By the above bound, $P_n^k[S]$ has an independent set $I$ of at least $|S|/(t+1)$ vertices. Now $I$ is also an independent set of $P_n^k$. Thus $|I|\leq n/(k+1)$. Hence $|S|/(t+1)\leq n/(k+1)$, and $|S|\leq (t+1)n/(k+1)$. \end{proof}

%%%%%%%%%%%%%%%%%%%%%%%%%%%%%%%%%%%%%%%%%%%%%%%%%%%%%%%%%%%%%%%%%%%%%%%%%%%%
\mySection{Structure of Bounded Degree Subgraphs}{Structure}
%%%%%%%%%%%%%%%%%%%%%%%%%%%%%%%%%%%%%%%%%%%%%%%%%%%%%%%%%%%%%%%%%%%%%%%%%%%%

In this section we study the structure of the subgraph of a $k$-tree induced by
the vertices of bounded degree. We first prove that in a $k$-tree with
sufficiently many vertices, not all the vertices of a clique have low degree. A
clique $C=(v_1,v_2,\dots,v_k)$ of a graph $G$ is said to be \emph{ordered by
degree} if $\deg_G(v_i)\leq\deg_G(v_{i+1})$ for all $1\leq i\leq k-1$. 

\begin{theorem}
\thmlabel{CliqueDegrees}
Let $G$ be a $k$-tree on $n\geq 2k+1$ vertices. Let $(u_1,u_2,\dots,u_q)$ be a
clique of $G$ ordered by degree. Then $\deg_G(u_i)\geq k+i-1$ for all $1\leq
i\leq q$. 
\end{theorem}

Note that \thmref{CliqueDegrees} is not true if $n\leq 2k$, as the statement
would imply that a ($k+1$)-clique has a vertex of degree $n$. Thus the
difficulty in an inductive prove of \thmref{CliqueDegrees} is the base case.
\thmref{CliqueDegrees} follows from the following stronger result with
$n\geq2k+1\geq k+q$.

\begin{lemma}
\lemlabel{VidasLemma}
Let $G$ be a $k$-tree on $n$ vertices.
Let $C=(u_1,u_2,\dots,u_q)$ be a clique of $G$ ordered by degree. 
If $n\geq k+q$ then
\begin{equation}
\eqnlabel{VidasLemmaA}
\deg_G(u_i)\,\geq\,k+i-1,\;\;1\leq i\leq q\enspace;
\end{equation}
otherwise $n\leq k+q-1$, and
\begin{equation}
\eqnlabel{VidasLemmaB}
\deg_G(u_i)\geq 
\begin{cases}
k+i-1	&\textup{ if }1\leq i\leq n-k-1\enspace, \\
n-1	&\textup{ if }n-k\leq i\leq q\enspace.
\end{cases}
\end{equation}
\end{lemma}

%%%%%%%%%%%%%%%%%%%%%

\begin{proof} We proceed by induction on $n$. In the base case, $G$ is
a $(k+1)$-clique, and every vertex has degree $k$. The claim follows trivially. 
Assume the result holds for $k$-trees on less than $n$ vertices. 
Let $C$ be a $q$-clique of a $k$-tree $G$ on $n\geq k+2$ vertices.
Since every $k$-tree on at least $k+2$
vertices has two non-adjacent simplicial vertices \citep{Dirac61}, 
at least one simplicial vertex $v$ is not in $C$. 
Since $n\geq k+2$ and $v$ is simplicial, the
graph $G_1=G\setminus v$ is a $k$-tree on $n-1$ vertices. 
Now $C$ is a $q$-clique of $G_1$. 
Let $C=(u_1,u_2,\dots,u_q)$ be ordered by degree in $G_1$.
By induction, if $n\geq k+q+1$ then
\begin{equation}
\eqnlabel{first}
\deg_{G_1}(u_i)\geq k+i-1,\; 1\leq i\leq q\;;
\end{equation}
otherwise $n\leq k+q$, and
\begin{equation}
\eqnlabel{second}
\deg_{G_1}(u_i)\geq
\begin{cases}
k+i-1	& \text{ if }1\leq i\leq n-k-2\enspace,\\ 
n-2		& \text{ if }n-k-1\leq i\leq q \enspace.
\end{cases}
\end{equation}

First suppose that $n\geq k+q+1$. Then by \eqnref{first},
$\deg_G(u_i)\geq\deg_{G_1}(u_i)\geq k+i-1$, and \eqnref{VidasLemmaA} is
satisfied. Otherwise $n\leq k+q$. Let $B=\{u_{n-k-1},u_{n-k},\dots,u_q\}$.
Then $|B|\geq 2$, and by \eqnref{second}, every vertex in $B$ has degree $n-2$
in $G_1$. That is, each vertex in $B$ is adjacent to every other vertex in
$G_1$. Let $X$ be the set of neighbours of $v$. Since $v$ is simplicial, $X$ is
a $k$-clique. At most one vertex of $B$ is not in $X$, as otherwise $X\cup B$ 
would be a ($k+2$)-clique of $G_1$. Without loss of generality, this
exceptional vertex in $B$, if it exists, is $u_{n-k-1}$. The other vertices in
$B$ are adjacent to one more vertex, namely $v$, in $G$ than in $G_1$. Thus 
$\deg_G(u_i)\geq k+i-1$ for all $1\leq i\leq n-k-1 $, and $\deg_G(u_i)=n-1$
for all $n-k\leq i\leq q$. Hence \eqnref{VidasLemmaB} is satisfied.
\end{proof} 

%%%%%%%%%%%%%%%%%%%%%%%%%%%%%%%%%%%%%%%%%%%%%%%%%%%%
We can now prove the main result of this section.

\begin{theorem}
\thmlabel{BoundedTreewidth}
For all integers $1\leq k\leq\ell\leq2k$, 
and for every $k$-tree $G$ on $n\geq\ell+2$ vertices, 
the subgraph $G_\ell$ of $G$ induced by the vertices of 
degree at most $\ell$, has treewidth at most $\ell-k$.
\end{theorem}

\begin{proof} Let $C=(u_1,u_2,\dots,u_q)$ be a clique of $G$ ordered by degree. Suppose, for the sake of contradiction, that there are at least $\ell-k+2$ vertices of $C$ with degree at most $\ell$. Let $j=\ell-k+2$. Since $C$ is ordered by degree, $\deg(u_j)\leq\ell$. Since $n\geq\ell+2$, we have $j\leq n-k$. By \lemref{VidasLemma}, $\deg(u_j)\geq k+j-1$ (unless $j=n-k$, in which case $\deg(u_j)=n-1\geq\ell+1$, which is a contradiction). Hence $k+j-1\leq\ell$. That is, $k+(\ell-k+2)-1\leq\ell$, a contradiction. Thus $C$ contributes at most $\ell-k+1$ vertices to $G_\ell$, and $\omega(G_\ell)\leq\ell-k+1$. Now, $G_\ell$ is an induced subgraph of $G$, which is chordal. Thus $G_\ell$ is chordal. Since $\omega(G_\ell)\leq\ell-k+1$, $G_\ell$ has treewidth at most $\ell-k$. \end{proof}

Note the following regarding \thmref{BoundedTreewidth}:

\begin{itemize}

\item There are graphs of treewidth $k\geq 2$ for which the theorem is not true. For example, for any $p\geq k+1$, consider the graph $G$ consisting of a ($k+1$)-clique $C$ and a $p$-vertex path with one endpoint $v$ in $C$. Then $G$ has at least $2k+1$ vertices, has treewidth $k$, and every vertex of $G$ has degree at most $k$, except for $v$ which has $\deg(v)=k+1$. For $\ell=k$, $G_\ell$ is comprised of two components, one a $k$-clique and the other a path, in which case $G_\ell$ has treewidth $k-1>\ell-k=0$. For $k+1\leq\ell\leq2k-1$, $G_\ell=G$ has treewidth $k>\ell-k$.

\item The theorem is not true if $k\leq n\leq\ell+1$. For example, for any
$1\leq k\leq \ell\leq 2k-1$, the $k$-tree obtained by adding $\ell+1-k$ vertices onto an initial $k$-clique has $\ell+1$ vertices, maximum degree $\ell$, and treewidth $k>\ell-k$.

\item The case of $\ell=k$ is the well-known fact that in a $k$-tree with at
least $k+2$ vertices, distinct simplicial vertices are not adjacent. Put
another way, the set of simplicial vertices of a $k$-tree with at least $k+2$
vertices is a $0$-set.

\end{itemize}

%%%%%%%%%%%%%%%%%%%%%%%%%%%%%%%%%%%%%%%%%%%%%%%%%%%%%%%%%%%%%%%%%
\mySection{Large Subgraphs of Bounded Treewidth and Bounded Degree}{LargeBoundedTreewidthDegreeSubgraphs}
%%%%%%%%%%%%%%%%%%%%%%%%%%%%%%%%%%%%%%%%%%%%%%%%%%%%%%%%%%%%%%%%%

The following theorem is the main result of the paper. 

\begin{theorem}
\thmlabel{dtSetLowerBound}
For all integers $0\leq t\leq k$, $d\geq2k$, and $n\geq 2k+1$,
\begin{equation*}
\alpha_d^t(\G_{n,k})\;\geq
\bracket{\frac{d-2k+1}{d-\threehalves k+1+\frac{t(t+1)}{2(k+1)}}}
\bracket{\frac{t+1}{k+1}}n +
 \frac{k(t+1)}{d-\threehalves k+2 + \frac{t(t+1)}{2(k+1)}}
\end{equation*}
\end{theorem}

\begin{proof}
Let $G$ be a graph in $\G_{n,k}$ with $\alpha_d^t(G)=\alpha_d^t(\G_{n,k})$. A
degree-$d$ $t$-set of a spanning supergraph of $G$ is a degree-$d$ $t$-set of
$G$. Thus we can assume that $G$ is a $k$-tree.

Consider $\ell$ with $k+t\leq\ell\leq 2k$. By \thmref{BoundedTreewidth}, $G_\ell$ has treewidth at most $\ell-k$. Since $t\leq\ell-k$, by \thmref{tSet},
\begin{equation*}
\alpha^t(G_\ell)\;\geq\;
\bracket{\frac{t+1}{\ell-k+1}}|V_\ell(G)|\enspace.
\end{equation*}
Since $\ell\leq d$, $\alpha^t(G_\ell)\leq\alpha_d^t(G)$, which implies that
\begin{equation}
\eqnlabel{Gen1}
|V_\ell(G)|\;\leq\;\bracket{\frac{\ell-k+1}{t+1}}\alpha_d^t(G)
\enspace.
\end{equation}

Now, $G$ has $kn-\half k(k+1)$ edges and minimum degree $k$. Let $n_i$ be the number of vertices of $G$ with degree exactly $i$. Thus,
\begin{equation*}
\sum_{i\geq k}i\cdot n_i
\;=\;2|E(G)|\;=\;
2kn-k(k+1)\;=\;-k(k+1)+\sum_{i\geq k}2k\cdot n_i\enspace.
\end{equation*}
Thus,
\begin{equation*}
\sum_{i\geq 2k+1}(i-2k)n_i
\;=\;-k(k+1)+\sum_{i=k}^{2k-1}(2k-i)n_i
\;=\;-k(k+1)+\sum_{i=k}^{2k-1}|V_i(G)|
\enspace,
\end{equation*}
and
\begin{equation*}
\sum_{i\geq 2k+1}(i-2k)n_i
\;=\;-k(k+1)+
\sum_{i=k}^{k+t-1}|V_i(G)|\;+\;\sum_{i=k+t}^{2k-1}|V_i(G)|\enspace.
\end{equation*}
By \eqnref{Gen1},

\begin{align*}
\sum_{i\geq 2k+1}(i-2k)n_i
&\;\leq\;-k(k+1)\;+\;
t\cdot |V_{k+t}(G)|\;+
\sum_{i=k+t}^{2k-1}\frac{(i-k+1)\cdot \alpha_d^t(G)}{t+1}\\
%%%%%%%%%%%%%%%%%%%%%%%%%%%%%
&\;\leq\;-k(k+1)\;+\;
t\cdot \alpha_d^t(G)\;+\;
\frac{\alpha_d^t(G)}{t+1}\sum_{i=t+1}^{k}i \\
%%%%%%%%%%%%%%%%%%%%%%%%%%%%%
&\;=\;-k(k+1)\;+\;
\alpha_d^t(G)\bracket{t\;+\;
\frac{1}{t+1}\bracket{\frac{k(k+1)-t(t+1)}{2}}}\\
%%%%%%%%%%%%%%%%%%%%%%%%%%%%%
&\;=\;-k(k+1)\;+\;\alpha_d^t(G)\bracket{\frac{t(t+1)+k(k+1)}{2(t+1)}}\enspace.
\end{align*}
Since $d\geq 2k$,
\begin{equation*}
-k(k+1)\;+\;\alpha_d^t(G)\bracket{\frac{t(t+1)+k(k+1)}{2(t+1)}}
\;\geq\;
\sum_{i\geq d+1}\!\!\!(i-2k)n_i
\;\geq\;
(d-2k+1)\sum_{i\geq d+1}\!\!\!n_i
\enspace.
\end{equation*}
Hence,
\begin{equation*}
|V_d(G)|\;=\;
n-\sum_{i\geq d+1}\!\!n_i
\;\geq\;
n\;+\;\frac{k(k+1)}{d-2k+1}
\;-\;\alpha_d^t(G)\bracket{\frac{t(t+1)+k(k+1)}{2(t+1)(d-2k+1)}}
\enspace.
\end{equation*}
By \thmref{tSet},
\begin{equation*}
\alpha_d^t(G)\;=\;
\alpha^t(G_d)\;\geq\;
\frac{t+1}{k+1}|V_{d}(G)|\;\geq\;
\frac{(t+1)n}{k+1}
\;+\;\frac{k(t+1)}{d-2k+1}
\;-\;\alpha_d^t(G)\bracket{\frac{t(t+1)+k(k+1)}{2(k+1)(d-2k+1)}}
\enspace.
\end{equation*}
Thus 
\begin{align*}
\bracket{1+\frac{t(t+1)+k(k+1)}{2(k+1)(d-2k+1)}}\alpha_d^t(G)
\;&\geq\;\frac{(t+1)n}{k+1}
\;+\;\frac{k(t+1)}{d-2k+1}\\
\bracket{\frac{2(k+1)(d-2k+1)+t(t+1)+k(k+1)}{2(k+1)(d-2k+1)}}\alpha_d^t(G)
\;&\geq\;\frac{(d-2k+1)(t+1)n\,+\,k(k+1)(t+1)}{(k+1)(d-2k+1)}\\
\Big((2d-3k+2)(k+1)+t(t+1)\Big)\alpha_d^t(G)
\;&\geq\;2(d-2k+1)(t+1)n\,+\,2k(k+1)(t+1)\\
\alpha_d^t(G)
\;&\geq\;\frac{2(d-2k+1)(t+1)n\,+\,2k(k+1)(t+1)}{(2d-3k+2)(k+1)+t(t+1)}\\
\alpha_d^t(G)
\;&\geq\;\frac{(d-2k+1)(t+1)n\,+\,k(k+1)(t+1)}
{(d-\threehalves k+1)(k+1)+\half t(t+1)}\enspace.
\end{align*}
The result follows.
\end{proof}

A number of notes regarding \thmref{dtSetLowerBound} are in order:

\begin{itemize}

\item \thmref{dtSetLowerBound} with $t=k$ is equivalent to the lower bound in \thmref{kSet}.

\item For $d<2k$, no result like \thmref{dtSetLowerBound} is possible, since $\alpha_d^t(P_n^k)=2(t+1)$.

\item The proof of \thmref{dtSetLowerBound} is similar to a strategy developed by \citet{BiedlWilkinson-ISAAC02-TCS} for finding bounded degree independent sets in planar graphs. 

\end{itemize}

%%%%%%%%%%%%%%%%%%%%%%%%%%%%%%%%%%%%%

\thmref{dtSetLowerBound} implies that there is a degree-$d$ $t$-set whose cardinality is arbitrarily close to the best possible bound without any degree restriction (\thmref{tSet}).

\begin{corollary}
For every $\epsilon>0$ and for all integers $0\leq t\leq k$, there exists
$d=d(\epsilon,k,t)$ such that for all $n\geq 2k+1$,
\begin{equation*}
\alpha_d^t(\G_{n,k})\;\geq\;(1-\epsilon)\bracket{\frac{t+1}{k+1}}n\enspace.
\end{equation*}
\end{corollary}

\begin{proof} By \thmref{dtSetLowerBound} it suffices to solve
\begin{equation*}
\frac{1-\epsilon}{k+1}
\;=\;
\frac{d-2k+1}{(d-\threehalves k+1)(k+1)+\half t(t+1)}\enspace.
\end{equation*}
That is,
%\begin{equation*}
%(1-\epsilon)\Big((d-\threehalves k+1)(k+1)+\half t(t+1)\Big)
%\;=\;
%(d-2k+1)(k+1)\enspace.
%\end{equation*}
%That is,
%\begin{equation*}
%-\epsilon(2d)(k+1)
%\;=\;
%(\epsilon-1)\big((-3k+2)(k+1)\,+t(t+1)\big)
%\;+(-4k+2)(k+1)
%\enspace,
%\end{equation*}
%and
\begin{equation*}
d\;=\;
\frac{1}{2}\bracket{1-\frac{1}{\epsilon}}
\bracket{3k-2-\frac{t(t+1)}{k+1}}
+\frac{2k-1}{\epsilon}
\enspace.
\end{equation*}
\end{proof}

%%%%%%%%%%%%%%%%%%%%%%%%%%%%%%%%%%%%%%%%%%%%%%%%%%%%%%%%%%%%%%%%%

We now prove an existential upper bound on the cardinality of a degree-$d$ $t$-set.

\begin{theorem}
\thmlabel{dtSetUpperBound}
For all integers $k\geq1$ and $d\geq 2k-1$ such that $2(d-2k+1)\equiv0\pmod{k(k+1)}$, there are infinitely many values of $n$, such that for all $0\leq t<k$,
\begin{equation*}
\alpha_d^t(\G_{n,k})\;\leq\;
\bracket{\frac{d-2k+1}{d-\threehalves k+1}}\bracket{\frac{t+1}{k+1}}n
\;+\;\frac{(k-1)(t+1)(d-2k+1)+k(t+1)(k+1)}{(d-\threehalves k+1)(k+1)}\enspace.
\end{equation*}
\end{theorem}

\begin{proof} Our construction employs the following operation. Let $G$ be a
$k$-tree containing an ordered $k$-clique $C=(v_1,v_2,\dots,v_k)$. A \emph{block} at $C$ consists of $k+1$ new vertices $\{x_1,x_2,\dots,x_{k+1}\}$ where $x_1$ is added onto the $k$-clique $\{v_1,v_2,\dots,v_k\}$; $x_2$ is added onto the $k$-clique $\{v_1,v_2,\dots,v_{k-1},x_1\}$; $x_3$ is added onto the $k$-clique $\{v_1,v_2,\dots,v_{k-2},x_1,x_2\}$; and so on, up to $x_{k+1}$ which is added onto the $k$-clique $\{x_1,x_2,\dots,x_k\}$. Clearly the graph obtained by adding a block to a $k$-clique of a $k$-tree is also a $k$-tree

Our graph is parameterised by the positive integer $n_0\geq 2k+3$. Initially let $G$ be the $k$-th power of a path $(v_1,v_2,\dots,v_{n_0})$. Note that any $k+1$ consecutive vertices in the path form a clique. Let $r$ be the non-negative integer such that $2(d-2k+1)=rk(k+1)$. Add $r$ blocks to $G$ at $(v_i,v_{i+1},\dots,v_{i+k-1})$ for each $3\leq i\leq n_0-k-1$, as illustrated in \figref{kTreeExample}.

\Figure[!htb]{kTreeExample}{\includegraphics{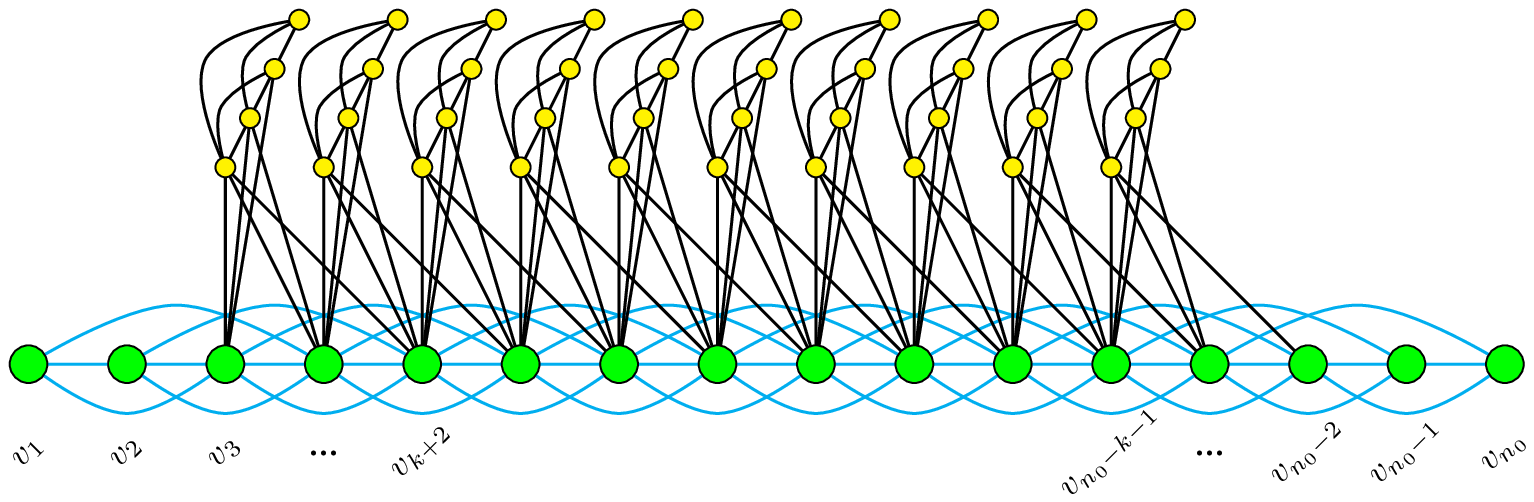}}{The graph $G$ with
$k=3$ and $d=11$ (and thus $r=1$).}

$G$ is a $k$-tree with $n=n_0+r(k+1)(n_0-(k+3))$ vertices. Let $S$ be a maximum degree-$d$ $t$-set of $G$. Consider a vertex $v_i$ for $k+2\leq i\leq n_0-k-1$. Since $n_0\geq 2k+3$ there is such a vertex. The degree of $v_i$ is
\begin{equation*}
2k+r\sum_{i=1}^ki=2k+\half rk(k+1)=d+1
\enspace.
\end{equation*}
Thus $v_i\not\in S$. Since each block $\{x_1,x_2,\dots,x_{k+1}\}$ is a clique, and treewidth-$t$ graphs have no $(t+2)$-clique, at most $t+1$ vertices from each block are in $S$. Similarly, since $\{v_1,v_2,\dots,v_{k+1}\}$ and $\{v_{n_0-k},v_{n_0-k+1},\dots,v_{n_0}\}$ are cliques, at most $t+1$ vertices from each of these sets are in $S$. Thus
\begin{equation}
\eqnlabel{ad}
\alpha_d^t(\G_{n,k})\;\leq\;
\alpha_d^t(G)\;=\;
|S|\;\leq\;
(t+1)\big(r(n_0-(k+3))+2\big)
\enspace.
\end{equation}
%The vertices $v_1$ and $v_{n_0}$ along with the vertex $x_{k+1}$ from each
%block form an independent set. Since $\deg(v_1)=k<d$, $\deg(v_{n_0})=k<d$, and
%$\deg(x_{k+1})=k<d$, we have equality in \eqnref{ad}. 
Substituting the equality $n_0=\frac{n+r(k+1)(k+3)}{1+r(k+1)}$
%\eqnref{n} 
into \eqnref{ad},
\begin{align}
\frac{\alpha_d^t(\G_{n,k})}{t+1}
%&\;\leq\;r(n_0-(k+3))+2\\
%&\;\leq\;r\Big(\frac{n+r(k+1)(k+3)}{1+r(k+1)}-(k+3)\Big)+2\\
%&\;\leq\;r\Big(\frac{n+r(k+1)(k+3) - (1+r(k+1))(k+3) }{1+r(k+1)}\Big)+2\\
%&\;\leq\;r\Big(\frac{n -(k+3) }{1+r(k+1)}\Big)+2\\
%&\;\leq\;\frac{r(n -(k+3)) +2 +2r(k+1)}{1+ r(k+1)}\\
%&\;\leq\;\frac{rn -rk -3r +2 +2rk +2r}{1+ r(k+1)}\\
%&\;\leq\;\frac{rn +rk -r +2 }{1+r(k+1)}\\
%&
\;\leq\;\frac{r(n+k-1) +2}{1+r(k+1)}\enspace.\eqnlabel{ada}
\end{align}
The claimed bound on $\alpha_d^t(\G_{n,k})$ follows by substituting the equality $r=\frac{2(d-2k+1)}{k(k+1)}$ into \eqnref{ada}.
%\begin{align*}
%\frac{\alpha_d^t(\G_{n,k})}{t+1}
%&\;\leq\;\frac{2\big(\frac{d-2k+1}{k(k+1)}\big)(n+k-1)+2}
%{1+2\big(\frac{d-2k+1}{k(k+1)}\big)(k+1)}\\
%&\;\leq\;\frac{\frac{2d-4k+2}{k(k+1)}(n+k-1) +2}{1+\frac{2d-4k+2}{k}}\\
%&\;\leq\;\frac{\frac{2d-4k+2}{k(k+1)}(n+k-1) +2}{\frac{2d-3k+2}{k}}\\
%&\;\leq\;\frac{(2d-4k+2)(n+k-1) +2k(k+1)}{(k+1)(2d-3k+2)}\\
%&\;\leq\;\frac{(2d-4k+2)n+(k-1)(2d-4k+2)+2k(k+1)}{(k+1)(2d-3k+2)}\enspace.
%\end{align*}
Observe that $n$ is a function of $n_0$ and $n_0$ is independent of $d$. Thus there are infinitely many values of $n$ for each value of $d$.
\end{proof}

%%%%%%%%%%%%%%%%%%%%%%%%%%%%%%%%%%%%%%%%%%%%%%%%%%%%%%%%%%%%%%%%%
\mySection{Bounded Degree Independent Sets}{IndSets}
%%%%%%%%%%%%%%%%%%%%%%%%%%%%%%%%%%%%%%%%%%%%%%%%%%%%%%%%%%%%%%%%%

Intuitively, one would expect that a maximum independent set would not have vertices $v$ of high degree, as this would prevent the many neighbours of $v$ from being in the independent set. In this section, we explore the accuracy of this intuition in the case of $k$-trees. Recall that $\alpha_d(G)$ is the maximum cardinality of a degree-$d$ independent set of $G$. 

%For a family of graphs \F,
%let $\alpha_d(\F)$ be the minimum, taken over all graphs $G\in\F$, of
%$\alpha_d(G)$. 
%Let $\PP_n$ be the family of $n$-vertex planar graphs. 
%
%that $\alpha_d(\PP_n)\in\Theta(1)$ for $d\leq 5$;
%$\alpha_d(\PP_n)=\frac{d-7}{4d-22}n+\Theta(1)$ for $6\leq d\leq 15$;
%$\frac{5}{23}n\leq \alpha_d(\PP_n)\leq \frac{d-7}{4d-22}n+2$ for $16\leq d\leq
%17$; $\frac{5}{23}n\leq \alpha_d(\PP_n)\leq \frac{3d-25}{12d-85}n+2$ for
%$d=18$; and $\frac{d-7}{4d-20}n\leq \alpha_d(\PP_n)\leq
%\frac{3d-25}{12d-85}n+2$ for $d\geq 19$.

Motivated by applications in computational geometry, the previously known results regarding bounded degree independent sets have been for planar graphs \citep{Kirkpatrick83, Edelsbrunner88, SvK-ESA97, BiedlWilkinson-ISAAC02-TCS}. The best results were obtained by \citet{BiedlWilkinson-ISAAC02-TCS}, who proved tight bounds (up to an additive constant) on $\alpha_d(G)$ for planar $G$ with $d\leq 15$. For $d\geq 16$ there is a gap in the bounds.

\thmref{tSet} with $t=0$ proves that every $n$-vertex graph $G$ with treewidth $k$ has $\alpha(G)\geq n/(k+1)$, and that this bound is tight for $P_n^k$. 
\thmref{dtSetLowerBound} with $t=0$ gives the following lower bound on the size of a degree-$d$ independent set in a graph of treewidth $k$ (for all $k\geq1$ and $d\geq2k$):
\begin{equation*}
\alpha_d(\G_{n,k})
\;\geq\;
\bracket{\frac{d-2k+1}{d-\threehalves k+1}}\bracket{\frac{n}{k+1}}\;+\;
\frac{k}{d-\threehalves k+1}\enspace.
\end{equation*}
Note that such a bound is not possible for $d<2k$ since $\alpha_d(P_n^k)=2$ for $d<2k$. 

\thmref{dtSetUpperBound} proves the corresponding upper bound. In particular, for all $k\geq1$, there are infinitely many values of $d$, and for each such $d$, there are infinitely many values of $n$ for which
\begin{equation*}
\alpha_d(\G_{n,k})\;\leq\;
\bracket{\frac{d-2k+1}{d-\threehalves k+1}}\bracket{\frac{n}{k+1}}
\;+\;\frac{(k-1)(d-2k+1)+k(k+1)}{(d-\threehalves k+1)(k+1)}\enspace.
\end{equation*}
These lower and upper bounds are tight. In fact, they differ by at most one. We conclude that 
\begin{equation*}
\lim_{n\rightarrow\infty}
\lim_{d\rightarrow\infty}
\frac{\alpha_d(\G_{n,k})}{n}
\;=\;
\frac{d-2k+1}{(d-\threehalves k+1)(k+1)}\enspace.
\end{equation*}

%%%%%%%%%%%%%%%%%%%%%%%%%%%%%%%%%%%%%%%%%%%%%%%%%%%%%%%%%%%%%%%%%%%%%%%%%%%%%%%
\subsection{Trees}
\seclabel{Trees}
%%%%%%%%%%%%%%%%%%%%%%%%%%%%%%%%%%%%%%%%%%%%%%%%%%%%%%%%%%%%%%%%%%%%%%%%%%%%%%%

$\G_{n,1}$ is precisely the family of $n$-vertex forests. Observe that \twothmref{dtSetLowerBound}{dtSetUpperBound} with $k=1$ and $t=0$ prove that for all $d\geq1$,
\begin{equation}
\eqnlabel{Trees}
\alpha_d(\G_{n,1})=\frac{(d-1)n+2}{2d-1}\enspace.
\end{equation}
A tree $T$ for which $\alpha_d(T)=\frac{(d-1)n+2}{2d-1}$ is called \emph{$\alpha_d$-extremal}. In this section we characterise the $\alpha_d$-extremal trees. A tree is \emph{$d$-regular} if every vertex has degree $1$ or $d$, and there is at least one vertex of degree $d$.

\begin{theorem}
\thmlabel{ExtremalTrees}
Let $d$ be a positive integer. A tree $T$ on $n\geq5$ vertices is $\alpha_d$-extremal if and only if $T$ is obtained from a $(d+1)$-regular tree by subdividing every leaf-edge once.
\end{theorem}

Note that a tree $T$ is $\alpha_1$-extremal if $\alpha_1(T)=2$. Every tree that is not a path has three independent leaves. Thus the only $\alpha_1$-extremal trees are paths, and \thmref{ExtremalTrees} holds trivially for $d=1$.
In the remainder of this section we consider the $d\geq 2$ case. We will use the following notation. For all trees $T$, let $L(T)$ be the set of leaves in $T$; let $P(T)$ be the set of degree-2 vertices in $T$; and let $Q(T)$ be the set of vertices in $P(T)$ that are not adjacent to a leaf. The following lemma is well-known.

\begin{lemma}
\lemlabel{RegularTreeLeaves}
For $d\geq2$, every $(d+1)$-regular tree $T$ with $n$ vertices satisfies \begin{equation*}
|L(T)|\;=\;\frac{(d-1)n+2}{d}\enspace.
\end{equation*}
\end{lemma}

\begin{proof} $T$ has $n-|L(T)|$ vertices of degree $d+1$ and has $n-1$ edges. Thus $|L(T)|+(d+1)(n-|L(T)|)=2(n-1)$. The result follows.
\end{proof}

\begin{proof}[Proof of \thmref{ExtremalTrees} \paran{$\Leftarrow$}] Let $T$ be a $(d+1)$-regular tree on $n$ vertices. By \lemref{RegularTreeLeaves}, $d|L(T)|=(d-1)n+2$. Let $T'$ be the tree obtained by subdividing every leaf-edge of $T$. Then $T'$ has $n'=n+|L(T)|$ vertices. Thus $d|L(T)|=(d-1)(n'-|L(T)|)+2$, which implies that $|L(T)|(2d-1)=(d-1)n'+2$. Now $T'$ has $2|L(T)|$ vertices of degree at most $d$, and they induce a matching. Thus $\alpha_d(T')=|L(T)|=((d-1)n'+2)/(2d-1)$, as claimed.
\end{proof}

We now prove a lower bound on $\alpha_d(T)$ that is more precise than  \thmref{dtSetLowerBound} with $k=1$ and $t=0$.  

\begin{lemma}
\lemlabel{IndSetTrees}
Let $T$ be a tree with $n\geq3$ vertices. Let $n_i$ be the number of vertices of $T$ with degree exactly $i$. For all $d\geq1$,
\begin{align*}
\alpha_d(T)\;\geq\;
\frac{1}{2d-1}\bracket{
(d-1)n+2
\;+\;\sum_{i=3}^{d}(i-2)n_i
\;+\;\sum_{i\geq d+2}(i-d-1)n_i}
\enspace.
\end{align*}
\end{lemma}

\begin{proof}
We proceed as in \thmref{dtSetLowerBound}. We have
\begin{equation*}
\sum_{i\geq1}i\cdot n_i
\;=\;2|E(T)|\;=\;
2n-2\;=\;-2+\sum_{i\geq 1}2n_i\enspace.
\end{equation*}
Thus,
\begin{align*}
n_1\;=\;2+\sum_{i\geq3}(i-2)n_i\enspace.
\end{align*}
Since $n\geq3$, no two leaves are adjacent. Thus $\alpha_d(T)\geq n_1$, and
\begin{align*}
\alpha_d(T)
\;&\geq\;
2+\sum_{i\geq3}(i-2)n_i\\
%%%%%%%%%%%%%%
%\alpha_d(T)
\;&=\;
2
\;+\;\sum_{i=3}^{d}(i-2)n_i
\;+\;\sum_{i\geq d+1}(i-d-1)n_i
\;+\;(d-1)\sum_{i\geq d+1}n_i\\
%%%%%%%%%%%%%%
%\alpha_d(T)
\;&=\;
2
\;+\;\sum_{i=3}^{d}(i-2)n_i
\;+\;\sum_{i\geq d+2}(i-d-1)n_i
\;+\;(d-1)(n-|V_d(T)|)\\
%%%%%%%%%%
(d-1)|V_d(T)|
\;&\geq\;
2-\alpha_d(T)
\;+\;(d-1)n
\;+\;\sum_{i=3}^{d}(i-2)n_i
\;+\;\sum_{i\geq d+2}(i-d-1)n_i
\enspace.
\end{align*}
The subgraph of $T$ induced by $V_d(T)$ is 2-colourable. 
The larger colour class is a degree-$d$ independent set of $T$. Thus
$\alpha_d(T)\geq\half|V_d(T)|$, which implies that
\begin{align*}
2(d-1)\alpha_d(T)
\;&\geq\;
2-\alpha_d(T)
\;+\;(d-1)n
\;+\;\sum_{i=3}^{d}(i-2)n_i
\;+\;\sum_{i\geq d+2}(i-d-1)n_i
\\
%%%%%%%%%%%%%%
(2d-1)\alpha_d(T)
\;&\geq\;
2
\;+\;(d-1)n
\;+\;\sum_{i=3}^{d}(i-2)n_i
\;+\;\sum_{i\geq d+2}(i-d-1)n_i
\enspace.
\end{align*}
The result follows.
\end{proof}

We have the following immediate corollary of \lemref{IndSetTrees}.

\begin{corollary}
\corlabel{ExtremalTreeDegrees}
For $d\geq1$, every vertex in an $\alpha_d$-extremal tree has degree in $\{1,2,d+1\}$.
\end{corollary}

%%%%%%%%%%%%%%%%%%%%%%%%%%%%%%%%%%%%%%%%%%%%%%%%%%%%%%%%%%%%%%%%%%%%%%%%%%%%%

\begin{lemma}
\lemlabel{ExtremalTreeInduction}
For $d\geq2$, every $n$-vertex tree $T$ in which every vertex has degree in $\{1,2,d+1\}$ satisfies
\begin{align*}
\alpha_d(T)\;&\geq\;
\frac{(d-1)n+\half|Q(T)|+2}{2d-1}\enspace.
\end{align*}
In particular, if $T$ is $\alpha_d$-extremal, then $Q(T)=\emptyset$.
\end{lemma}

\begin{proof}
By \lemref{RegularTreeLeaves} applied to the $(d+1)$-regular tree obtained from $T$ by contracting every vertex of degree two, 
\begin{equation*}
|L(T)|\;=\;\frac{(d-1)(n-|P(T)|)+2}{d}\enspace.
\end{equation*}
There is at most one vertex in $P(T)$ adjacent to each leaf. Thus $|P(T)|\leq|L(T)|+|Q(T)|$. Hence
\begin{equation*}
|L(T)|\;\geq\;\frac{(d-1)(n-|L(T)|-|Q(T)|)+2}{d}\enspace,
\end{equation*}
which implies that
%d|L(T)|\;&\geq\;(d-1)(n-|L(T)|-|Q(T)|)+2\\
%(2d-1)|L(T)|\;&\geq\;(d-1)n-(d-1)|Q(T)|+2\\
\begin{equation*}
|L(T)|\;\geq\;\frac{(d-1)n-(d-1)|Q(T)|+2}{2d-1}\enspace.
\end{equation*}
The subgraph of $T$ induced by $Q(T)$ is a forest of paths, no vertex of which is adjacent to a leaf. Thus $\alpha_d(T)\geq|L(T)|+\half|Q(T)|$. Hence,
\begin{equation*}
\alpha_d(T)
\;\geq\;\frac{(d-1)n-(d-1)|Q(T)|+2}{2d-1}+\frac{|Q(T)|}{2}
\;=\;\frac{(d-1)n+\half|Q(T)|+2}{2d-1}
\enspace,
\end{equation*}
as desired.
\end{proof}

%%%%%%%%%%%%%%%%%%%%%%%%%%%%%%%%%%%%%%%%%%%%%%%%

\begin{lemma}
\lemlabel{ExtremalTreePT}
For $d\geq2$, every $\alpha_d$-extremal tree $T$ on $n\geq5$ vertices has $|P(T)|\geq|L(T)|$.
\end{lemma}

\begin{proof}
By \corref{ExtremalTreeDegrees}, every vertex in $T$ has degree in $\{1,2,d+1\}$. $T$ is not a path as otherwise $\alpha_d(T)\geq\half n>((d-1)n+2)/(2d-1)$ for $n\geq5$. Let $T'$ be the tree obtained from $T$ by contracting every vertex of degree two. Then $T'$ is $(d+1)$-regular, and has $n'=n-|P(T)|$ vertices and $|L(T)|$ leaves. By \lemref{RegularTreeLeaves},
\begin{equation}
\eqnlabel{Sizen1}
|L(T)|=\frac{(d-1)n'+2}{d}\enspace.
\end{equation}
Clearly, $\alpha_d(T)\geq |L(T)|$. Since $T$ is $\alpha_d$-extremal,
\begin{equation*}
\frac{(d-1)n'+2}{d}\;=\;|L(T)|
\;\leq\;\alpha_d(T)\;=\;\frac{(d-1)n+2}{2d-1}\enspace.
\end{equation*}
Now $n=n'+|P(T)|$. Thus
\begin{align*}
\frac{(d-1)n'+2}{d}\;&\leq\;\frac{(d-1)(n'+|P(T)|)+2}{2d-1}\\
(2d-1)(d-1)n'+2(2d-1)\;&\leq\;d(d-1)(n'+|P(T)|)+2d\\
(d-1)^2n'+2(d-1)\;&\leq\;d(d-1)|P(T)|\\
(d-1)n'+2\;&\leq\;d|P(T)|\enspace.
\end{align*}
By \eqnref{Sizen1}, $|P(T)|\geq|L(T)|$, as claimed.
\end{proof}

%%%%%%%%%%%%%%%%%%%%%%%%%%%%%%%%%%%%%%%%%%%

\begin{proof}[Proof of \thmref{ExtremalTrees} \paran{$\Rightarrow$}] Let $T$ be an $\alpha_d$-extremal tree. By \corref{ExtremalTreeDegrees}, every vertex of $T$ has degree in $\{1,2,d+1\}$. By \lemref{ExtremalTreeInduction}, $Q(T)=\emptyset$. That is, every degree-2 vertex is adjacent to a leaf. By \lemref{ExtremalTreePT}, $|P(T)|\geq|L(T)|$. That is, there are at least as many degree-2 vertices as leaves. Hence $|P(T)|=|L(T)|$, and $T$ is obtained from a $(d+1)$-regular tree by subdividing every leaf-edge once.
\end{proof}

%\COMMENT{Is there a function $f(d)<1$ such that for every tree $T$ we have
%$\alpha_d(T)\geq f(d)\alpha(T)$ and $f(d)\rightarrow1$ for large $d$? Use %dynamic programming.}

%%%%%%%%%%%%%%%%%%%%%%%%%%%%%%%%%%%%%%%%%%%%%%%%%%%%%%%%%%%%%%%%%%%%%%%%%%%%%%%
\subsection{Outerplanar Graphs}
\seclabel{Outerplanar}
%%%%%%%%%%%%%%%%%%%%%%%%%%%%%%%%%%%%%%%%%%%%%%%%%%%%%%%%%%%%%%%%%%%%%%%%%%%%%%%

A plane embedding of a graph in which every vertex is on a single face is called \emph{outerplanar}. A graph is \emph{outerplanar} if it has an outerplanar embedding. Let $\OP_n$ denote the class of $n$-vertex outerplanar graphs. It is well known that the outerplanar graphs are a proper subset of the class of graphs with treewidth at most two (see \citep{Bodlaender-TCS98}). However, the graphs constructed in the upper bound in \thmref{dtSetUpperBound} with $k=2$ are not outerplanar. We have the following upper bound for outerplanar graphs.

\begin{theorem}
\thmlabel{Outerplanar}
For all $d\geq4$ and $n\geq5$,
\begin{equation*}
\alpha_d(\OP_n)\;\geq
\bracket{\frac{d-3}{3d-6}}n +  \frac{2}{d-2}\enspace.
\end{equation*}
Conversely, for all even $d\geq6$ and for infinitely many values of $n$,
\begin{equation*}
\alpha_d(\OP_n)\;\leq
\bracket{\frac{d-4}{3d-10}}\bracket{n-6}+3\enspace.
\end{equation*}
\end{theorem}

\begin{proof}
The lower bound follows from \thmref{dtSetLowerBound} with $k=2$.

For the upper bound, let $r=(d-4)/2$. Since $d\geq6$ is even, $r$ is a positive integer. Our graph $G$ is parameterised by an integer $n_0\geq6$. Initially let $G=P^2_{n_0}$ be the square of a path $(v_1,v_2,\dots,v_{n_0})$. That is, $v_iv_j$ is an edge whenever $|i-j|\leq2$. For each vertex $v_i$, $3\leq i\leq n_0-2$, $G$ has $r$ triangles $\{a_{i,j},b_{i,j},c_{i,j}\}$, $1\leq j\leq r$, where $v_i$ is adjacent to each $a_{i,j}$ and $b_{i,j}$. In addition, for $3\leq i\leq n_0-4$, the edge $v_ia_{i+2,1}$ is in $G$. Finally there are two additional vertices $x$ and $y$; $x$ is adjacent to $v_{n_0-3}$ and $v_{n_0-1}$, and $y$ is adjacent to $v_{n_0-2}$ and $v_{n_0}$. $G$ has $n=n_0+(n_0-4)3r+2$ vertices. As illustrated in \figref{OuterplanarExample}, there is an outerplanar embedding of $G$.

\Figure[!htb]{OuterplanarExample}{\includegraphics{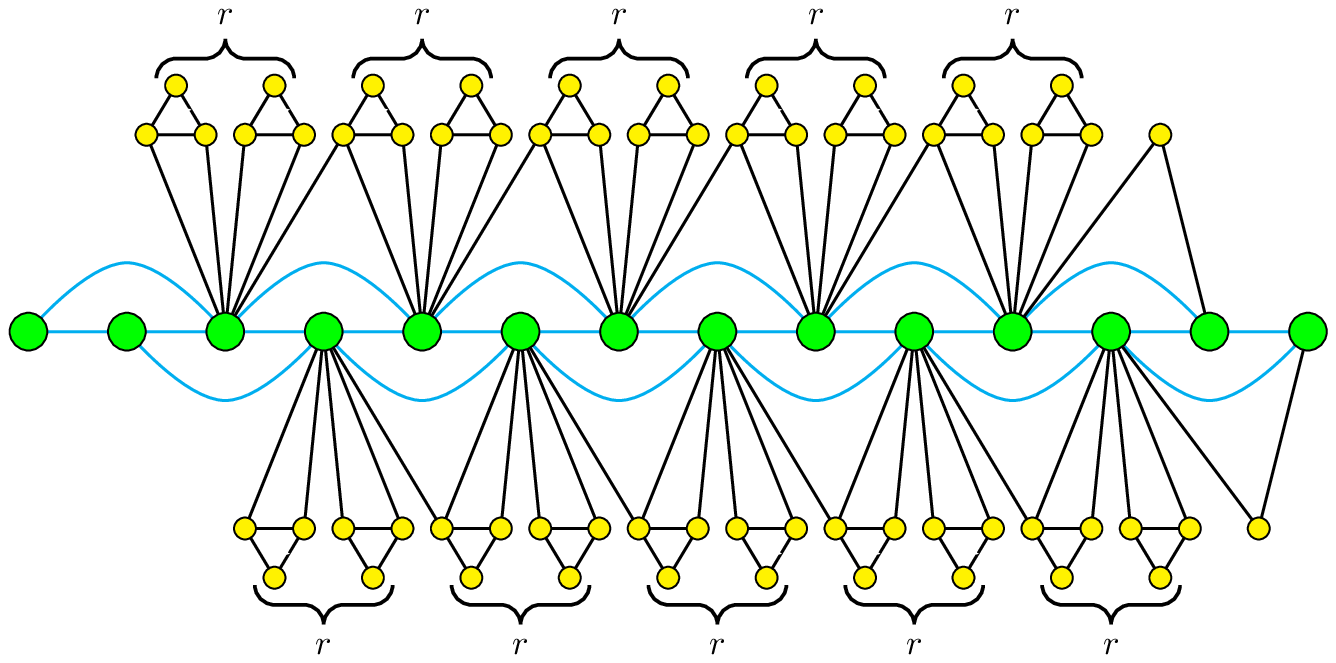}}{The outerplanar graph $G$ for $r=2$.}

%\COMMENT{Do you see a nicer way to do this? What about odd $d$? I can do %$\frac{d-3}{3d-13}n$ for odd $d$.}

Say $S$ is a degree-$d$ independent set of $G$. Each vertex $v_i$, $3\leq i\leq n_0-2$, has degree $4+2r+1=d+1$, and is thus not in $S$. At most one vertex from each triangle $\{a_{i,j},b_{i,j},c_{i,j}\}$ is in $S$. It follows that $|S|\leq r(n_0-4)+3$. Now $$n_0=\frac{n+12r-2}{3r+1}.$$ 
Thus 
\begin{equation*}
|S|
\;\leq\;
r\bracket{\frac{n+12r-2}{3r+1}-4}+3
%&=\;
%r\bracket{\frac{n-6}{3r+1}}+3\\
%&=\;
%\frac{r}{3r+1}\bracket{n-6}+3\\
%&=\;
%\frac{(d-4)/2}{3(d-4)/2+1}\bracket{n-6}+3\\
%&
\;=\;
\frac{d-4}{3d-10}\bracket{n-6}+3\enspace.
\end{equation*}
The result follows.
\end{proof}

Note that the upper and lower bound in \thmref{Outerplanar} are tight for
$d=6$. That is, every $n$-vertex outerplanar graph $G$ has an degree-$6$ independent  set on $\quarter n -\Oh{1}$ vertices, and for infinitely many values of $n$, there is an $n$-vertex outerplanar graph in which at most $\quarter n+\Oh{1}$ vertices form a degree-$6$ independent set. Recall that without any degree restriction, every outerplanar graph has an independent set on at least $\third n$ vertices. An interesting open problem is to derive  upper and lower bounds on $\alpha_d(\OP_n)$ that are tight for infinitely  many values of $d$.
 
%%%%%%%%%%%%%%%%%%%%%%%%%%%%%%%%%%%%%%%%%%%%%%%%%%%%%%%%%%%%%%%%%%%%%%%%%%%%%%%
\subsection{Interval Graphs}
\seclabel{Pathwidth}
%%%%%%%%%%%%%%%%%%%%%%%%%%%%%%%%%%%%%%%%%%%%%%%%%%%%%%%%%%%%%%%%%%%%%%%%%%%%%%%

A graph $G$ is an \emph{interval graph} if one can assign to each vertex $v\in
V(G)$ a closed interval $[L_v,R_v]\subseteq\mathbb{R}$ such that $vw\in E(G)$
if and only if $[L_v,R_v]\cap[L_w,R_w]\ne\emptyset$. An interval graph $G$ has tree-width equal to $\omega(G)+1$. (In fact, it has path-width equal to $\omega(G)+1$.) Thus the previous results of this paper apply to interval graphs. However, for bounded degree independent sets in interval graphs, we can say much more, as we show in this section. In an interval graph, it is well known that we can assume that the endpoints of the intervals are distinct. We say a vertex $w$ is \emph{dominated} by a vertex $v$ if $L(v)<L(w)<R(w)<R(v)$. 

\begin{lemma}
\lemlabel{CoveredVertex}
Let $G$ be an interval graph with $\omega(G)\leq k+1$. Suppose $G$ has a vertex $v$ with $\deg(v)\geq2k+1$. Then there is a vertex $w$ that is dominated by $v$ and $\deg(w)\leq 2k-1$.
\end{lemma}

\begin{proof} 
For each vertex $y\in V(G)$, let $A(y)=\{x\in V(G):L(x)<L(y)<R(x)\}$ and $B(y)=\{x\in V(G):L(x)<R(y)<R(x)\}$. Observe that $x$ is dominated by $y$ if and only if $xy\in E(G)$ but $x\not\in A(y)\cup B(y)$. Also $|A(y)|\leq k$ as otherwise $A(y)\cup\{y\}$ would be a clique of at least $k+2$ vertices. Similarly $|B(y)|\leq k$. Thus $|A(y)\cup B(y)|\leq 2k$. 

Now consider the given vertex $v$. Since $\deg(v)\geq 2k+1$, $v$ has a neighbour $u\not\in A(v)\cup B(v)$. Thus $u$ is dominated by $v$. Let $w$ be a vertex with the shortest interval that is dominated by $v$. That is, if $u$ and $w$ are dominated by $v$, then $R(w)-L(w)\leq R(u)-L(u)$. Thus $w$ does not dominate any vertex, and every neighbour of $w$ is in $A(w)\cup B(w)$. Now $|A(w)|\leq k$,  $|B(w)|\leq k$, and $v\in A(w)\cap B(w)$. Thus $\deg(w)\leq2k-1$.  \end{proof}

Note that \lemref{CoveredVertex} with $k=1$ is the obvious statement that a vertex of degree at least three in a caterpillar is adjacent to a leaf.

\begin{theorem}
\thmlabel{IntervalGraphs}
Every interval graph $G$ with $\omega(G)\leq k+1$ has a degree-$2k$ maximum
independent set. That is, $\alpha_{2k}(G)=\alpha(G)$.
\end{theorem}

\begin{proof} Let $I$ be a maximum independent set of $G$. If $I$ contains a
vertex $v$ with $\deg(v)\geq 2k+1$, apply \lemref{CoveredVertex} to obtain a
vertex $w$ dominated by $v$ such that $\deg(w)\leq 2k-1$. Replace $v$ by $w$ in
$I$. The obtained set is still independent, since every neighbour of $w$ is
also adjacent to $v$, and is thus not in $I$. Apply this step repeatedly until
every vertex in $I$ has degree at most $2k$. Thus
$\alpha_{2k}(G)\geq|I|=\alpha(G)$. By definition, $\alpha_{2k}(G)\leq
\alpha(G)$. Therefore $\alpha_{2k}(G)=\alpha(G)$. \end{proof}

The bound of $2k$ in \thmref{IntervalGraphs} is best possible, since $P_n^k$ is
an interval graph with $\omega(G)\leq k+1$ and only $2k$ vertices of degree at most $2k-1$. Thus $\alpha(P_n^k)=\ceil{n/(k+1)}\gg\alpha_{2k-1}(P_n^k)$.

%%%%%%%%%%%%%%%%%%%%%%%%%%%%%%%%%%%%%%%%%%%%%%%%%%%%%%%%%%%%%%%%%%%%%%%%%%%%%%%
%\bibliographystyle{myNatbibStyle}
%\bibliography{myBibliography,myConferences}

\begin{thebibliography}{8}
\providecommand{\natexlab}[1]{#1}
\expandafter\ifx\csname urlstyle\endcsname\relax
  \providecommand{\doi}[1]{doi:\discretionary{}{}{}#1}\else
  \providecommand{\doi}{doi:\discretionary{}{}{}\begingroup
  \urlstyle{rm}\Url}\fi

\bibitem[{Biedl and Wilkinson(2002)}]{BiedlWilkinson-ISAAC02-TCS}
\textsc{Therese Biedl and Dana~F. Wilkinson}.
\newblock Bounded-degree independent sets in planar graphs.
\newblock In \textsc{Prosenjit Bose and Pat Morin}, eds., \emph{Proc. 13th
  International Conf. on Algorithms and Computation (ISAAC '02)}, vol. 2518 of
  \emph{Lecture Notes in Comput. Sci.}, pp. 416--427. Springer, 2002.
\newblock To appear in \emph{Theory Comput. Syst.}

\bibitem[{Bodlaender(1998)}]{Bodlaender-TCS98}
\textsc{Hans~L. Bodlaender}.
\newblock A partial $k$-arboretum of graphs with bounded treewidth.
\newblock \emph{Theoret. Comput. Sci.}, 209(1-2):1--45, 1998.

\bibitem[{Bose \emph{et~al.}(2004)Bose, Smid, and Wood}]{BSW-DM04}
\textsc{Prosenjit Bose, Michiel Smid, and David~R. Wood}.
\newblock Light edges in degree-constrained graphs.
\newblock \emph{Discrete Math.}, 282(1-3):35--41, 2004.

\bibitem[{Dirac(1961)}]{Dirac61}
\textsc{Gabriel~A. Dirac}.
\newblock On rigid circuit graphs.
\newblock \emph{Abh. Math. Sem. Univ. Hamburg}, 25:71--76, 1961.

\bibitem[{Edelsbrunner(1988)}]{Edelsbrunner88}
\textsc{Herbert Edelsbrunner}.
\newblock \emph{Algorithms in Combinatorial Geometry}.
\newblock Springer, 1988.

\bibitem[{Kirkpatrick(1983)}]{Kirkpatrick83}
\textsc{David Kirkpatrick}.
\newblock Optimal search in planar subdivisions.
\newblock \emph{SIAM J. Comput.}, 12(1):28--35, 1983.

\bibitem[{Reed(2003)}]{Reed-AlgoTreeWidth03}
\textsc{Bruce~A. Reed}.
\newblock Algorithmic aspects of tree width.
\newblock In \textsc{Bruce~A. Reed and Cl{\'a}udia~L. Sales}, eds.,
  \emph{Recent Advances in Algorithms and Combinatorics}, pp. 85--107.
  Springer, 2003.

\bibitem[{Snoeyink and van Kreveld(1997)}]{SvK-ESA97}
\textsc{Jack Snoeyink and Marc van Kreveld}.
\newblock Linear-time reconstruction of {D}elaunay triangulations with
  applications.
\newblock In \textsc{Rainer~E. Burkhard and Gerhard~J. Woeginger}, eds.,
  \emph{Proc. 5th Annual European Symp. on Algorithms (ESA '97)}, vol. 1284 of
  \emph{Lecture Notes in Comput. Sci.}, pp. 459--471. Springer, 1997.

\end{thebibliography}
%%%%%%%%%%%%%%%%%%%%%%%%%%%%%%%%%%%%%%%%%%%%%%%%%%%%%%%%%%%%%%%%%%%%%%%%%%%%%%%%

\def\soft#1{\leavevmode\setbox0=\hbox{h}\dimen7=\ht0\advance \dimen7
  by-1ex\relax\if t#1\relax\rlap{\raise.6\dimen7
  \hbox{\kern.3ex\char'47}}#1\relax\else\if T#1\relax
  \rlap{\raise.5\dimen7\hbox{\kern1.3ex\char'47}}#1\relax \else\if
  d#1\relax\rlap{\raise.5\dimen7\hbox{\kern.9ex \char'47}}#1\relax\else\if
  D#1\relax\rlap{\raise.5\dimen7 \hbox{\kern1.4ex\char'47}}#1\relax\else\if
  l#1\relax \rlap{\raise.5\dimen7\hbox{\kern.4ex\char'47}}#1\relax \else\if
  L#1\relax\rlap{\raise.5\dimen7\hbox{\kern.7ex
  \char'47}}#1\relax\else\message{accent \string\soft \space #1 not
  defined!}#1\relax\fi\fi\fi\fi\fi\fi} \def\cprime{$'$}

\end{document}